\documentclass[12pt]{amsart}
\usepackage{amssymb}
\usepackage{fullpage}

\newcommand{\g}{\mathfrak g}
\renewcommand{\t}{\mathfrak t}
\renewcommand{\b}{\mathfrak b}
\newcommand{\n}{\mathfrak n}
\renewcommand{\u}{\mathfrak u}

\newcommand{\s}{\mathrm s}

\renewcommand{\P}{\mathbb{P}}

\newcommand{\cB}{\mathcal B}
\newcommand{\cC}{\mathcal C}
\newcommand{\cN}{\mathcal N}
\newcommand{\cO}{\mathcal O}

\DeclareMathOperator{\PGL}{PGL} %
\DeclareMathOperator{\SL}{SL} %
\DeclareMathOperator{\Sp}{Sp} %
\DeclareMathOperator{\SO}{SO} %

\DeclareMathOperator{\Char}{char} %
\DeclareMathOperator{\Lie}{Lie} %
\DeclareMathOperator{\rank}{rank} %

\newtheorem{thm}{Theorem}[section]
\newtheorem{lem}[thm]{Lemma}
\newtheorem{cor}[thm]{Corollary}

\theoremstyle{definition}
\newtheorem{exmp}[thm]{Example}
\newtheorem{prob}[thm]{Problem}

\theoremstyle{remark}
\newtheorem{rem}[thm]{Remark}

\title[The orbit structure of Dynkin curves]
{The orbit structure of Dynkin curves}
\author[S.~M.~Goodwin, L.~Hille and G.~R\"ohrle]
{Simon M.~Goodwin, Lutz Hille and Gerhard R\"ohrle}

\address{School of Mathematics, University of Birmingham,
Birmingham, B15 2TT, UK} \email{goodwin@maths.bham.ac.uk}
\urladdr{http://web.mat.bham.ac.uk/S.M.Goodwin/}

\address{Mathematisches Seminar,
Universit\"at Hamburg, D-20 146, Hamburg, Germany}
\email{hille@math.uni-hamburg.de}
\urladdr{http://www.math.uni-hamburg.de/home/hille/}

\address{School of Mathematics, University of Southampton,
Southampton, SO17 1BJ, UK} \email{G.Roehrle@soton.ac.uk}
\urladdr{http://www.maths.soton.ac.uk/staff/Roehrle}


\dedicatory{To the memory of Peter Slodowy}

\begin{document}

\begin{abstract}
Let $G$ be a simple algebraic group over an algebraically closed
field $k$; assume that $\Char k$ is zero or good for $G$.  Let $\cB$
be the variety of Borel subgroups of $G$ and let $e \in \Lie G$ be
nilpotent. There is a natural action of the centralizer $C_G(e)$ of
$e$ in $G$ on the Springer fibre $\cB_e = \{B' \in \cB \mid e \in
\Lie B'\}$ associated to $e$. In this paper we consider the case,
where $e$ lies in the subregular nilpotent orbit; in this case
$\cB_e$ is a Dynkin curve.  We give a complete description of the
$C_G(e)$-orbits in $\cB_e$. In particular, we classify the
irreducible components of $\cB_e$ on which $C_G(e)$ acts with
finitely many orbits.
In an application we obtain a classification of all subregular
orbital varieties admitting a finite number of $B$-orbits for $B$ a
fixed Borel subgroup of $G$.
\end{abstract}

\maketitle

\section{Introduction}

Let $G$ be a simple algebraic group over the algebraically closed
field $k$; assume that the characteristic of $k$ is zero or good for
$G$. Let $e \in \g = \Lie G$ be nilpotent and write $\cO = G \cdot
e$ for the adjoint $G$-orbit of $e$. Let $B$ be a Borel subgroup of
$G$ and let $\u$ be the Lie algebra of the unipotent radical of $B$;
we may assume that $e \in \u$ . We write $\cB$ for the variety of
Borel subgroups in $G$ and let $\cB_e = \{ B' \in \cB \mid e \in
\u'\}$ be the Springer fibre associated to $e$.  We note that
$\cB_e$ is connected (\cite[Prop.\ 1]{St}),
equidimensional (\cite[Prop.\ 1.12]{Sp}) and of dimension
$\frac{1}{2}(\dim C_G(e) - \rank G)$ (\cite[\S 5]{Ca}).

The varieties $\cB_e$ occur as the fibres in Springer's resolution
$\{(e',B') \in \cN \times \cB \mid e' \in \u'\} \to \cN$
of the nilpotent variety $\cN$ of $\g$.  This desingularization,
which is the projection map onto $\cN$, was given by T.~A.~Springer
in \cite{Spr}, see also \cite[6.10]{Ja}. The Springer fibres have
attracted much research interest due to their importance in the
representation theory of the Weyl group $W$ of $G$; more
specifically, Springer constructed all irreducible representations
of $W$ on the top cohomology groups of $\cB_e$, see \cite{Spr2} and
\cite{Lu}. Subsequently there has been much interest in the
geometric structure of the varieties $\cB_e$, see for example
\cite{DLP}, \cite{HS} and \cite{St2}.

Since all the Borel subgroups of $G$ are conjugate, we have $\cB =
\{{}^gB \mid g \in G\}$ is the variety of all $G$-conjugates of $B$.
There is a well-known connection between the Springer fibre $\cB_e$
and the variety $\cO \cap \u$ due to N.~Spaltenstein, see \cite[II
1.9]{Sp}. Let $\pi_1 : G \to \cB$ be the map $g \mapsto
{}^{g^{-1}}B$ and let $\pi_2 : G \to \cO$ be the orbit map $g
\mapsto g\cdot e$. There is an action of $C_G(e) \times B$ on $Y =
\pi_1^{-1}(\cB_e) = \pi_2^{-1}(\cO \cap \u)$ by $(g,b) \cdot y =
byg^{-1}$. There a bijection between the $C_G(e)$-orbits in $\cB_e$
and the $(C_G(e) \times B)$-orbits in $Y$, and a bijection between
the $B$-orbits in $\cO \cap \u$ and the $(C_G(e) \times B)$-orbits
in $Y$. Therefore, there is a bijection between the $C_G(e)$-orbits
in $\cB_e$ and the $B$-orbits in $\cO \cap \u$; moreover this
bijection preserves the closure order on orbits. Further, the orbits
of the component group $A(e)=C_G(e)/C_G(e)^\circ$ of $C_G(e)$ on the
set of irreducible components of $\cB_e$ are in bijective
correspondence with the irreducible components of $\cO \cap \u$.

\smallskip

Due to the importance of the Springer fibres alluded to above it is
a natural problem to try to understand the orbits of $C_G(e)$ in
$\cB_e$ under the conjugation action.  This is formalized in the
problem posed below.

\begin{prob}
\label{Prob:fibre} Determine the orbit structure for the action of
$C_G(e)$ on $\cB_e$. In particular, determine when $C_G(e)^\circ$
acts on a given irreducible component of $\cB_e$ with a finite
number of orbits and if this is not the case, then determine whether
there is still a dense $C_G(e)^\circ$-orbit.
\end{prob}

It follows from the above discussion that Problem \ref{Prob:fibre}
is equivalent to Problem \ref{Prob:richardson} below.

\begin{prob}
\label{Prob:richardson} Determine the orbit structure for the action
of $B$ on $\cO \cap \u$.  In particular, determine when  $B$ acts on
an irreducible component of $\cO \cap \u$ with a finite number of
orbits and if this is not the case, then determine whether there is
still a dense $B$-orbit.
\end{prob}

First we record two elementary cases for Problem \ref{Prob:fibre}.
If $e$ is in the regular nilpotent orbit, then $\cB_e$ is a point
(cf.\ \cite[3.7 Thm.\ 1]{St}) and there is nothing to show. For $e = 0$,
we have $\cB_e = \cB$ is an irreducible homogeneous variety for
$C_G(e) = G$.

Thanks to the equivalence between these two problems we immediately
obtain a partial answer of Problem \ref{Prob:fibre} in the case of
spherical nilpotent orbits. Suppose that $\cO$ is a \emph{spherical}
nilpotent orbit; that is one on which $B$ admits a dense orbit.
Thanks to a fundamental result, due to M.~Brion \cite{Br} and
E.~B.~Vinberg \cite{Vi} in characteristic $0$ and to F.~Knop
\cite[2.6]{Kn} in arbitrary characteristic,  $B$ acts on a spherical
orbit $\cO$ with a finite number of orbits. Therefore, it is clear
that there are only finitely many $B$-orbits in each irreducible
component of $\cO \cap \u$. Correspondingly, in the spherical case
there are only finitely many $C_G(e)$-orbits in $\cB_e$.  All
spherical nilpotent orbits for $\Char k = 0$ have been classified by
D.~I.~Panyushev in \cite{Pa}. This classification has recently been shown to
hold for positive good characteristic in work of R.~Fowler and the
third author \cite{FR}.

\smallskip

The purpose of this paper is to give complete answers to Problems
\ref{Prob:fibre} and \ref{Prob:richardson} when $e$ lies in the
subregular nilpotent orbit, see Theorems \ref{T:fibre} and
\ref{T:richardson}.  When $G$ is simply laced, this classification
is nicely expressed in terms of coefficients of the highest root
$\rho$ of $G$.  As explained below, $\cB_e$ is a Dynkin curve. In
case $G$ is simply laced, the irreducible components of $\cB_e$ are
indexed by the simple roots of $G$, and Theorem \ref{T:fibre} says
that there is a finite number of $C_G(e)$-orbits in the irreducible
component corresponding to the simple root $\alpha$ if and only if
the coefficient of $\alpha$ in $\rho$ is 1. For $G$ non-simply laced,
the classification requires the notion of associated simply laced
root systems, see \S \ref{sub:curves}.

\smallskip

Let $\cO$ be a nilpotent $G$-orbit.  The irreducible components of
the variety $\overline{\cO \cap \u}$ are called \emph{orbital
varieties} and are of interest in the representation theory of the
Weyl group of $G$, see for example \cite[\S 9.13]{Ja}. They are also
of interest in the study of primitive ideals of the universal
enveloping algebra $U(\g)$ of $\g$, see for example \cite{BB} and
\cite{Jo}.  As an application of Theorem \ref{T:richardson}, we
classify when there is a finite number of $B$-orbits in an
irreducible component of $\overline{\cO \cap \u}$, when $\cO$ is the
subregular nilpotent orbit, see Theorem \ref{T:orbital}.

\smallskip

From now on let $e$ be a representative of the subregular nilpotent
orbit $\cO$ in $\g$. We recall from \cite[\S 3.10]{St}
that in this case $\cB_e$ is a {\em Dynkin curve}.
A Dynkin curve is a non-empty,
connected union of certain projective lines
determined by the root system of $G$. Each of the projective lines
has a type $\alpha$, where $\alpha$ is a simple root of $G$; they
are denoted by $\cC_\alpha^i$, where $i =
1,\dots,|\alpha|^2/|\alpha_\s|^2$ (here $\alpha_\s$ is a fixed short
root). We refer the reader to Section \ref{S:fibre} for a precise
definition of a Dynkin curve.

Dynkin curves are of interest as they arise in resolutions of
Kleinian singularities.  More specifically, a Dynkin curve occurs as
the exceptional divisor in the minimal resolution of a Kleinian
singularity.  We refer the reader to Slodowy's book \cite[\S 6]{Sl}
for details.

\smallskip

As a consequence of our solution to Problem \ref{Prob:fibre}, we
deduce that there are finitely many $C_G(e)$-orbits in $\cB_e$ if
and only if $G$ is of type $A_r$ or $B_r$; in both of these cases
there are exactly $2r-1$ orbits.  Further, we see that the action of
$C_G(e)$ on $\cB_e$ is trivial if and only if $G$ is of type $E_8$;
see Corollary \ref{cor:1}.

Since the irreducible components of $\cB_e$ are projective lines
$\cC_\alpha^i \cong \P^1$, if $C_G(e)^\circ$ acts on $\cC_\alpha^i$
with a dense orbit, then its complement in $\cC_\alpha^i$ is closed
and so it is a finite set of points. This implies that in the
subregular case, $C_G(e)^\circ$ admits a dense orbit in $\cC_\alpha^i$
precisely when there are finitely many orbits and so these two parts
of Problem \ref{Prob:fibre} are equivalent in this instance;
correspondingly, this is also the case for Problem
\ref{Prob:richardson}.

\smallskip

We remark that Problems \ref{Prob:fibre} and \ref{Prob:richardson}
could be stated for connected reductive $G$. We choose not to work
in this generality, as one can easily reduce to the case when $G$ is
simple, so there is no loss in generality. Further, we could
consider the analogous problems for unipotent $G$-conjugacy classes
in place of nilpotent $G$-orbits. Under our assumption that $\Char
k$ is zero or a good prime for $G$, there exists a \emph{Springer
isomorphism}, see \cite[III, 3.12]{SS} and \cite[Cor.\ 9.3.4]{BR},
so that these problems are equivalent.

\smallskip

For general references on reductive algebraic groups and nilpotent
orbits the reader is referred to Borel's book \cite{Bo} and
Jantzen's monograph \cite{Ja}.

\section{Preliminaries and statements of results}
\label{S:fibre}

\subsection{Notation}
Let $G$ be a simple algebraic group over an algebraically closed
field $k$; we assume throughout that the characteristic of $k$ is
zero or a good prime $p$ for $G$. Let $B$ be a Borel subgroup of $G$
with unipotent radical $U$ and let $T$ be a maximal torus of $G$
contained in $B$. We write $\g$, $\b$, $\u$ and $\t$ for the Lie
algebras of $G$, $B$, $U$ and $T$, respectively.  We write $\cB$ for
the variety of Borel subgroups of $G$.  Since all Borel subgroups of
$G$ are conjugate and $B$ is self-normalizing, we may identify $\cB$
with $G/B$.

Let $\Phi = \Phi(G,T)$ be the root system of $G$ with respect to
$T$. Let $\Phi^+$ be the system of positive roots determined by $B$
and $\Pi$ the corresponding base.  We denote the Dynkin diagram of
$\Phi$ by $\Delta$; we identify the nodes in $\Delta$ with the
elements of $\Pi$. Let $\rho$ be the highest root of $\Phi$ with
respect to $\Pi$. The Cartan matrix of $\Phi$ is written as
$(\langle \alpha, \beta \rangle)_{\alpha, \beta \in \Pi}$. Let
$\alpha_\s$ denote a fixed root of short length (if there is only
one root length in $\Phi$, then all roots are short). We label the
simple roots in the Dynkin diagram $\Delta$ in accordance with
\cite[Planches I--IX]{Bou}.

For a simple root $\alpha \in \Pi$, we let $e_\alpha$ be a generator
of the root subspace of $\g$ corresponding to $\alpha$. We write
$P_\alpha$ for the minimal parabolic subgroup of $G$ corresponding
to $\alpha$, i.e.\ $P_\alpha$ is the parabolic subgroup of $G$
generated by $B$ and the root subgroup
corresponding to $-\alpha$.  We write $\u_\alpha$ for
the Lie algebra of the unipotent radical of $P_\alpha$.

For the rest of this section we let $e$ be a representative of the
subregular nilpotent orbit $\cO = G \cdot e$; we assume that $e \in
\u$. We write $C_G(e)$ for the centralizer of $e$ in $G$ and $A(e) =
C_G(e)/C_G(e)^\circ$ for the component group of $C_G(e)$.  We denote
the Springer fibre associated to $e$ by $\cB_e = \{B' \in \cB \mid e
\in \u'\}$.

\subsection{Dynkin curves} \label{sub:curves}
We recall, from \cite[\S 3.10]{St} (see also \cite[\S 6.3]{Sl}),
that the Springer fibre $\cB_e$ is the \emph{Dynkin curve}
determined by the root system $\Phi$. This Dynkin curve $\cC$ is a
non-empty union of projective lines $\cC_\alpha^i$, where $\alpha
\in \Pi$ is the {\em type} of $\cC_\alpha^i$ and $i =
1,\dots,|\alpha|^2/|\alpha_\s|^2$. A line of type $\alpha$ meets
exactly $-\langle \alpha, \beta \rangle$ lines of type $\beta \ne
\alpha$ and they intersect in a single point; lines of the same type
do not meet.  We write $\cC_\alpha$ for $\cC_\alpha^1$;
$\cC_\alpha'$ for $\cC_\alpha^2$ and $\cC_\alpha''$ for
$\cC_\alpha^3$. Under the identification of $\cB$ with $G/B$ lines
of type $\alpha$ are of the form $xP_\alpha/B$, for some $x \in G$,
see \cite[\S 6.3]{Sl}.

It is straightforward to explicitly describe the Dynkin curve in
case $G$ is simply laced: there is a single projective line
$\cC_\alpha$ for each $\alpha \in \Pi$; and $\cC_\alpha$ intersects
$\cC_\beta$ if and only if $\alpha$ and $\beta$ are adjacent in
$\Delta$.  In other words, the Dynkin curve $\cC$ is the diagram
\emph{dual} to the Dynkin diagram $\Delta$, see \cite[\S 6.3]{Sl}.

In order to give an explicit description of $\cC$ in case $G$ is
non-simply laced we need to recall some further notation from
\cite[\S 6.2]{Sl}. Suppose $\Delta$ is non-simply laced, then we
define the {\em associated simply laced diagram} $\widehat \Delta$
of $\Delta$ and the \emph{associated symmetry group}
$\Gamma(\Delta)$ by Table~\ref{t:2} below; in this table $S_n$
denotes the symmetric group of degree $n$.
\begin{table}[h!]
\renewcommand{\arraystretch}{1.5}
\begin{tabular}{|c|cccc|}
\hline $\Delta$ & $B_r$ & $C_r$ & $F_4$ & $G_2$ \\
\hline $\widehat \Delta$ & $A_{2r-1}$ & $D_{r+1}$ & $E_6$ & $D_4$ \\
\hline $\Gamma(\Delta)$ & $S_2$ & $S_2$ & $S_2$ & $S_3$ \\ \hline
\end{tabular}
\vspace{2mm}
\caption{The associated Dynkin diagrams} \label{t:2}
\end{table}
\vspace{-7mm}

We write $\widehat \Phi$ for the root system and $\widehat \Pi$ for
the base of $\widehat \Phi$ corresponding to $\widehat \Delta$, and
$\widehat \rho$ for the highest root of $\widehat \Phi$ with respect
to $\widehat \Pi$. As explained in \cite[\S 6.2]{Sl}, there is a
unique faithful action of $\Gamma(\Delta)$ on $\widehat \Delta$ and
we can regard $\Delta$ as the quotient of $\widehat \Delta$ by this
action; we refer the reader to \cite[\S 6.2]{Sl} for a description
of the $\Gamma(\Delta)$-orbit in $\widehat \Pi$ corresponding to a
simple root in $\Pi$.

Suppose $G$ is non-simply laced. Then the Dynkin curve $\cC$ can be
explicitly described as follows.  Let $\widehat \cC$ be the Dynkin
curve determined by the simply laced root system $\widehat \Phi$. As
a variety $\cC$ is the same as $\widehat \cC$.  A projective line in
$\cC$ is of type $\alpha$ if its type in $\widehat \cC$ is in the
$\Gamma(\Delta)$-orbit corresponding to $\alpha$.

The action of $G$ on $\cB$ restricts to an action of $C_G(e)$ on
$\cB_e$.  This in turn induces an action of the component group
$A(e)$ of $C_G(e)$ on the set of lines of type $\alpha$ in $\cB_e$.
The group $\Gamma(\Delta)$ is precisely the component group $A(e)$
of $C_G(e)$; further, the action of $\Gamma(\Delta)$ on $\widehat
\Delta$ corresponds in a natural way to the action of $A(e)$ on the
lines in $\cC$, see \cite[Prop.\ 7.5]{Sl}. More precisely, there is
an isomorphism $\phi : \Gamma(\Delta) \to A(e)$ such that for $\beta
\in \widehat \Phi$ and $\cC_\beta$ the line of type $\beta$ in
$\widehat \cC$, and $g \in \Gamma(\Delta)$, we have $\phi(g) \cdot
\cC_{\beta} = \cC_{g \cdot \beta}$.  In particular, the action of
$A(e)$ on lines of type $\alpha$ is transitive.

We illustrate the above discussion with an example, cf.\
\cite[6.3]{Sl}.
\begin{exmp}
Consider the case when $G$ is of type $G_2$.
Let $\alpha = \alpha_1$ and $\beta = \alpha_2$ be the short and long
simple roots of $\Pi$, respectively. The Cartan matrix is
$\left(\begin{tabular}{rr} $2$ & $-1$
\\ $-3$ & $2$\end{tabular}\right)$. Therefore, $\cC$ is of the form
$\cC = \cC_\alpha \cup \cC_\beta \cup \cC_\beta' \cup \cC_\beta''$
and the intersection pattern is illustrated below. \unitlength 0.7pt

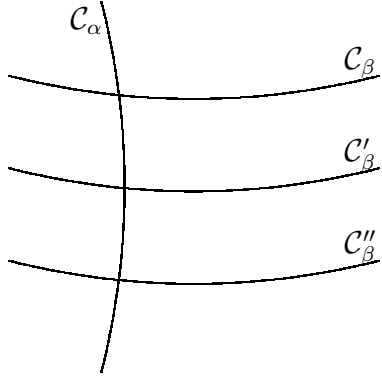
\[
\begin{picture}(200,200)
\qbezier(50,0)(75,100)(50,200)
\put(42,190){\makebox(0,0){{$\cC_\alpha$}}}
\qbezier(0,160)(100,135)(200,160)
\put(190,167){\makebox(0,0){{$\cC_\beta$}}}
\qbezier(0,110)(100,85)(200,110)
\put(190,117){\makebox(0,0){{$\cC_\beta'$}}}
\qbezier(0,60)(100,35)(200,60)
\put(190,67){\makebox(0,0){{$\cC_\beta''$}}}
\end{picture}
\]

As follows from the discussion above, as an algebraic variety,
the Dynkin curve for $G_2$ is the same as that for $C_3$ and $D_4$;
although the types of the lines differ.
\end{exmp}

For fixed $\alpha \in \Pi$, let $\Upsilon_\alpha$ be the partition
of $\cC_\alpha$ afforded by the points in the intersections
$\cC_\alpha \cap \cC_\beta$, $\cC_\alpha \cap \cC_\beta'$,
$\cC_\alpha \cap \cC_\beta''$  for all $\beta \in \Pi\setminus \{\alpha\}$,
and the
complement in $\cC_\alpha$ of the union of all these points;
$\Upsilon'_\alpha$ and $\Upsilon''_\alpha$ are defined analogously.
We observe that the partition of all of $\cB_e$ that we obtain by
taking the union of all the $\Upsilon_\alpha$, $\Upsilon'_\alpha$
and $\Upsilon''_\alpha$ for all $\alpha \in \Pi$
is related to the \emph{component configuration}
of $\cB_e$ as defined in \cite[1.4]{BH}.

We can now state the main theorem of this article; it is proved in
the next section.  For $\alpha \in \Pi$, we let $\widehat \alpha \in
\widehat \Pi$ be a representative of the $\Gamma(\Delta)$-orbit
corresponding to $\alpha$.

\begin{thm}
\label{T:fibre} Let $e$ be a representative of the subregular
nilpotent $G$-orbit in $\g$.  Then the following hold:
\begin{itemize}
\item[(a)]
There is a finite number of $C_G(e)^\circ$-orbits in $\cC_\alpha$ if
and only if
\begin{itemize}
\item[(i)]  the coefficient of $\alpha$ in $\rho$ is
$1$, for $G$ simply laced;
\item[(ii)] the coefficient of $\widehat \alpha$ in $\widehat
\rho$ is $1$, for $G$ non-simply laced.
\end{itemize}
\item[(b)]
If $C_G(e)^\circ$ acts on $\cC_\alpha$ with a finite number of
orbits, then each component of the partition $\Upsilon_\alpha$ of
$\cC_\alpha$ is a single $C_G(e)^\circ$-orbit.  Otherwise the action
of $C_G(e)^\circ$ on $\cC_\alpha$ is trivial.
\end{itemize}
\end{thm}

The following is immediate from Theorem \ref{T:fibre}.

\begin{cor}
\label{cor:1} Let $e$ be a representative of the subregular
nilpotent $G$-orbit in $\g$.  Then the following hold:
\begin{itemize}
\item[(i)]
There are finitely many $C_G(e)$-orbits in $\cB_e$ if and only if
$G$ is of type $A_r$ or $B_r$. Moreover, in both of these cases the
number of $C_G(e)$-orbits is $2r-1$.
\item[(ii)]
The action of $C_G(e)$ on $\cB_e$ is trivial if and only if $G$ is
of type $E_8$.
\end{itemize}
\end{cor}

\subsection{The varieties $\cO \cap \u$}
\label{sub:richardson}
We now discuss the varieties $\cO \cap \u$
for $\cO = G \cdot e$ the subregular nilpotent orbit.

The irreducible components of $\cO \cap \u$ are the intersections
$\cO \cap \u_\alpha$ for $\alpha \in \Pi$.  Under the correspondence
between the $C_G(e)$-orbits in $\cB_e$ and the $B$-orbits in $\cO
\cap \u$, the $C_G(e)$-orbits in lines of type $\alpha$ correspond
to the $B$-orbits in $\cO \cap \u_\alpha$, i.e.\ on the level of
correspondence of irreducible components: the $A(e)$-orbit of lines
of type $\alpha$ corresponds to $\cO \cap \u_\alpha$.  This allows
one to derive from the knowledge of the intersection pattern of
projective lines in the Dynkin curve $\cC$ that the intersection of
$\cO$ with $\u_\alpha \cap \u_\beta$ is either empty or a single
$B$-orbit. More precisely, $\cO \cap (\u_\alpha \cap \u_\beta) \ne
\varnothing$, if and only if $\langle \alpha, \beta \rangle \ne 0$,
and in this case this intersection is a single $B$-orbit.
This fact can also be deduced from
Richardson's dense orbit theorem (\cite{Ri}),
as explained in the next paragraph.

It is well-known that the subregular class $\cO$ is the
\emph{Richardson class} of $P_\alpha$ for each $\alpha \in \Pi$,
i.e.\ $\cO$ meets each $\u_\alpha$ in an open dense subvariety and
$\cO \cap \u_\alpha$ is a single $P_\alpha$-orbit. Observe that if
$\langle \alpha, \beta \rangle = 0$, then $\cO \cap (\u_\alpha \cap
\u_\beta) = \varnothing$, as in this case $\u_\alpha \cap \u_\beta$
is the Lie algebra of the unipotent radical of a larger parabolic
subgroup of $G$, which has a different Richardson class. It follows
from Lemma \ref{L:dimension} later in this paper that $\cO \cap
(\u_\alpha \cap \u_\beta)$ can be at most one $B$-orbit; and in case
$\langle \alpha, \beta \rangle \ne 0$ one can show $\cO \cap
(\u_\alpha \cap \u_\beta) \ne \varnothing$ by using the root
subgroups corresponding to $\beta$ and $-\beta$ to find a
representative of the Richardson orbit in $\u_\beta$ for which the
coefficient of $e_\alpha$ is zero.

Let $\Xi_\alpha$ be the partition of $\cO \cap \u_\alpha$ afforded
by all the subvarieties $\cO \cap (\u_\alpha \cap \u_\beta)$ for
all $\beta \in \Pi$ with $\langle \alpha, \beta \rangle \ne 0$ along
with the complement in $\cO \cap \u_\alpha$ of the union of all
these.

Now we state the counterpart to Theorem \ref{T:fibre} in the context
of the action of $B$ on $\cO \cap \u$; due to the discussion in the
introduction, these two theorems are equivalent.  In the next
section we prove both theorems by proving complementary parts of
each.

\begin{thm}
\label{T:richardson}
Let $\cO$ be the subregular nilpotent $G$-orbit
in $\g$. Then the following hold:
\begin{itemize}
\item[(a)]
There is a finite number of $B$-orbits in $\cO \cap \u_\alpha$ if
and only if
\begin{itemize}
\item[(i)]  the coefficient of $\alpha$ in $\rho$ is
$1$, for $G$ simply laced;
\item[(ii)] the coefficient of $\widehat \alpha$ in $\widehat
\rho$ is $1$, for $G$ non-simply laced.
\end{itemize}
\item[(b)]
If $B$ acts on $\cO \cap \u_\alpha$ with a finite number of orbits,
then each component of the partition $\Xi_\alpha$ of $\cO \cap
\u_\alpha$ is a single $B$-orbit.
\end{itemize}
\end{thm}

We may now state an analogue of Corollary \ref{cor:1}(i) in the
setting of $B$ acting on $\cO \cap \u$.

\begin{cor}
\label{cor:2}
Let $\cO$ be the subregular nilpotent $G$-orbit in
$\g$.  Then there are finitely many $B$-orbits in $\cO \cap \u$ if
and only if $G$ is of type $A_r$ or $B_r$. Moreover, in both of
these cases the number of $B$-orbits is $2r-1$.
\end{cor}

\section{Proofs of Theorems \ref{T:fibre} and \ref{T:richardson}}
\label{S:proof}

Before we begin the proof of Theorems \ref{T:fibre} and
\ref{T:richardson}, we state the following well-known lemma about
the action of tori and unipotent groups on projective lines.  It can
be deduced easily from \cite[Prop.\ 10.8]{Bo}, which says that any
action of an algebraic group $H$ on $\P^1$ is given by a
homomorphism $H \to \PGL_2(k)$.

\begin{lem} \label{L:action}
Let $H$ be a connected algebraic group acting non-trivially on
$\P^1$.
\begin{enumerate}
\item[(i)] Suppose $H$ is unipotent. Then $H$ has two orbits in $\P^1$:
one is a fixed point and the other is the complement of the fixed
point.
\item[(ii)]  Suppose $H$ is a torus.
Then $H$ has three orbits in $\P^1$: two are fixed points and the
other is their complement.
\end{enumerate}
\end{lem}

We recall that a nilpotent element $e \in \g$ is called
\emph{distinguished} if $C_G(e)^\circ$ is unipotent, see for example
\cite[\S 4.1]{Ja}.
Our next lemma follows directly from
\cite[\S 7.5 Lem.\ 4]{Sl}.

\begin{lem}
\label{L:distinguished} The subregular nilpotent class $\cO = G
\cdot e$ of $G$ is distinguished unless $G$ is of type $A_r$ or
$B_r$, in which case the Levi factor of $C_G(e)^\circ$ is a
one-dimensional torus.
\end{lem}

Armed with Lemmas \ref{L:action} and \ref{L:distinguished}, we can
deduce that the action of $C_G(e)^\circ$ on $\cC_\alpha$ ($\alpha
\in \Pi$) is trivial in many cases. These are the instances when $e$
is distinguished and $\cC_\alpha$ meets at least two other lines.
Since the connected centralizer $C_G(e)^\circ$ fixes these
intersection points and is unipotent, it has to act trivially on all
of $\cC_\alpha$, by Lemma \ref{L:action}(i).

This leaves us to consider the following cases: $G$ is of type $A_r$
and $B_r$ and $\alpha$ is any simple root;
or $G$ is of any type and $\alpha$ is an end-node of $\Delta$.

\smallskip

First, we use computations explained in \cite{Go} to complete the
proof of Theorem \ref{T:richardson} (equivalently Theorem
\ref{T:fibre}) for $G$ of exceptional type.
It follows from the results of these computations that $B$ acts on
$\cO \cap \u_\alpha$ with a dense orbit precisely in the cases
stated in Theorem \ref{T:richardson}. This in turn implies that
$C_G(e)^\circ$ admits a dense orbit in $\cC_\alpha$ in the stated
cases. Finally, we use the fact that $C_G(e)^\circ$ is unipotent and
Lemma \ref{L:action}(i) to deduce that $C_G(e)^\circ$ acts on
$\cC_\alpha$ with two orbits if it acts with a dense orbit. This
completes the discussion of the exceptional cases.

\smallskip

Next we complete the proof of Theorems \ref{T:fibre} and
\ref{T:richardson} for $G$ of classical type.  First, we explain why
we can reduce to considering groups of matrices.

Let $\phi: G \to \tilde G$ be an isogeny. Then it is clear that
Theorems \ref{T:fibre} and \ref{T:richardson} hold for $G$ if and
only if they hold for $\tilde G$.  Therefore,
we may assume $G$ is one of $\SL_n(k)$, $\SO_{2n+1}(k)$,
$\Sp_{2n}(k)$ or $\SO_{2n}(k)$. In these cases the nilpotent
$G$-orbits are parameterized by the partitions given by the Jordan
normal form (with some exceptions in case $G = \SO_{2n}(k)$ and $n$
is even, which are not relevant for our purposes here), see for
example \cite[Chapter 1]{Ja}. We require the following well-known
lemma, giving the partitions corresponding to subregular elements,
see for example \cite[4.5.6, Cor.\ 2]{He} and \cite[IV 2.33]{SS}.

\begin{lem} \label{L:subregular}
The partition corresponding to the subregular nilpotent class is:
\begin{enumerate}
\item[(i)]  $(n-1,1)$ for $G = \SL_n(k)$;
\item[(ii)] $(2n-1,1,1)$ for $G = \SO_{2n+1}(k)$;
\item[(iii)] $(2n-2,2)$ for $G = \Sp_{2n}(k)$;
\item[(iv)] $(2n-3,3)$ for $G = \SO_{2n}(k)$.
\end{enumerate}
\end{lem}

The following lemma, which holds for an arbitrary simple algebraic
group $G$, is used extensively to prove existence of a dense
$B$-orbit in $\cO \cap \u_\alpha$ in the classical cases.  It can be
proved using a simple dimension argument using the two equalities:
$\dim P_\alpha \cdot e = \dim \u_\alpha$ and $\dim B = \dim P_\alpha
- 1$; we omit the details.

\begin{lem} \label{L:dimension}
Let $\alpha \in \Pi$ and suppose $e \in \u_\alpha$ is subregular.
Then $\dim B \cdot e = \dim \u_\alpha$ or $\dim \u_\alpha - 1$.
Thus if $\dim C_B(e) < \dim C_{P_\alpha}(e)$, then $B \cdot e$ is
dense in $\u_\alpha$.
\end{lem}

We now prove that there is a dense $B$-orbit in $\u_\alpha$ ($\alpha
\in \Pi$) for the cases stated in Theorem \ref{T:richardson}.

\smallskip

First we consider the case $G = \SL_n(k)$.  We note that in this
case the existence of a dense $B$-orbit in $\u_\alpha$ can be
deduced directly from the main theorem in \cite{GH}, but we give a
more elementary proof here. We take $T$ to be the maximal torus of
diagonal matrices in $G$ and $B$ to be the Borel subgroup of upper
triangular matrices in $G$.
We write $e_{i,j}$ for the elementary matrix with $(i,j)$th entry
$1$ and all other entries $0$. Given a root $\alpha \in \Phi$, we let
$x_\alpha : k \to G$ be a parametrization of the corresponding root
subgroup of $G$ and we write $e_\alpha$ for a generator of the root
subspace, so $e_\alpha = e_{i,j}$ for some $i,j$.

Suppose $\alpha = \alpha_i$, where $i \ne 1, n-1$. Consider
\[
e =  \left(\sum_{j=1, j \ne i}^{n-1} e_{j,j+1}\right)  + e_{i,i+2} =
 \left(\sum_{j=1, j \ne i}^{n-1} e_{\alpha_j}\right)
+ e_{\alpha_i + \alpha_{i+1}}.
\]
One can check that $e^{n-2} \ne 0$ and $e^{n-1} = 0$. So, using
Lemma \ref{L:subregular}, we deduce that $e$ lies in the subregular
nilpotent orbit.  Consider $x_{-\alpha}(s) \cdot e$. One can see
that there exists $t_s \in T$ such that $t_s x_{-\alpha}(s) \cdot e
= e$ for all but one value of $s$. This implies that $\dim C_B(e) <
\dim C_{P_\alpha}(e)$. So, by Lemma \ref{L:dimension}, we have that
$B \cdot e$ is dense in $\u_\alpha$. From the equivalence of
Theorems \ref{T:fibre} and \ref{T:richardson}, we deduce there is a
dense $C_G(e)$-orbit in $\cC_\alpha$.  It now follows from Lemma
\ref{L:action} that there are precisely three $C_G(e)$-orbits.  This
completes the proof of Theorem \ref{T:fibre}, and therefore also
Theorem \ref{T:richardson}, in this case.

Now suppose $\alpha = \alpha_1$ (the case $\alpha = \alpha_{n-1}$ is
equivalent, by symmetry).  Consider
\[
e = \sum_{j=2}^{n-1} e_{j,j+1} = \sum_{j=2}^{n-1} e_{\alpha_j}
\]
and let $\cO = G \cdot e$. One can check that $e^{n-2} \ne 0$ and
$e^{n-1} = 0$.  So $e$ lies in the subregular nilpotent orbit. We
see that $x_{-\alpha}(s) \in C_G(e)$ for all $s \in k$, so it
follows from Lemma \ref{L:dimension} that $B \cdot e$ is dense in
$\u_\alpha$. Now, using the results of \cite[\S 7]{Go2}, one can show
that
$$
B \cdot e = \left\{x \in \u_\alpha \:\: \vline \:\: x = \sum_{i=2}^n
a_i e_{i,i+1} + \sum_{j=2}^n \sum_{i=1}^{j-2} b_{i,j} e_{i,j},\  a_i
\in k^\times, b_{i,j} \in k \right\}.
$$
Therefore, we have
$$
\u_\alpha \setminus (B \cdot e) = \bigcup \u_\alpha \cap \u_\beta,
$$
where the union is taken over all simple roots $\beta \ne \alpha$.

If $\beta \ne \alpha_2$, then $(G \cdot e) \cap (\u_\alpha \cap
\u_\beta) = \varnothing$, as in this case $\u_\alpha \cap \u_\beta$
is the Lie algebra of the unipotent radical of a parabolic subgroup
bigger than $P_\alpha$.  Using Lemma \ref{L:dimension}, we see that
$\cO \cap (\u_\alpha \cap \u_\beta)$ is a single $B$-orbit when
$\beta = \alpha_2$. So there are precisely two $B$-orbits in
$\cO \cap \u_\alpha$.
This completes the proof of Theorem \ref{T:richardson} in this case.

\smallskip

We now consider the other classical groups; the proofs in these
cases are similar to the one for $\SL_n(k)$.  We just give
representatives of the dense $B$-orbit in $\u_\alpha$ in the cases
stated in Theorem \ref{T:richardson}; one can prove that the
representatives do indeed give a dense orbit and that there is the
right number of $B$-orbits in $\cO \cap \u_\alpha$ using arguments
similar to those for the  $\SL_n(k)$-case, so we omit the
details.

Let $G$ be one of $\SO_{2n+1}(k)$, $\Sp_{2n}(k)$ or $\SO_{2n}(k)$.
Let $V$ be the natural $G$-module with standard
(ordered) basis $v_1,\dots,v_n,v_0,v_{-n},\dots,v_{-1}$ and
$G$-invariant symmetric or alternating bilinear form $(\,,)$ defined
by $(v_0, v_i) = (v_0, v_{-i}) = 0, (v_0,v_0) = 1$, $(v_i, v_j) =
(v_{-i},v_{-j}) = 0$ and $(v_i, v_{-j}) = \delta_{i,j}$ for $i,j =
1, \dots, n$, where we omit $v_0$ everywhere if $G$ is $\Sp_{2n}(k)$
or $\SO_{2n}(k)$.

We take $T$ to be the maximal torus of diagonal matrices in $G$ and
$B$ to be the Borel subgroup of upper triangular matrices in $G$.

\smallskip

First we consider $G = \SO_{2n+1}(k)$ and the case where
$\alpha = \alpha_1$. Then
\[
e =  \left(\sum_{j = 2}^{n-1} e_{j,j+1} - e_{-(j+1),-j}\right) +
(e_{n,0} - e_{0,-n}) =  \sum_{j = 2}^n e_{\alpha_j}
\]
is a representative of the dense $B$-orbit in $\u_\alpha$. For
$\alpha = \alpha_i$ for $i = 2,\dots,n-1$,
\begin{align*}
e = & \left(\sum_{j = 1, j \ne i}^{n-1} e_{j,j+1} -
e_{-(j+1),-j}\right) +
(e_{n,0} - e_{0,-n}) + (e_{i,i+2} - e_{-(i+2),-i}) \\
= & \left(\sum_{j = 1, j \ne i}^n e_{\alpha_j}\right) + e_{\alpha_i
+ \alpha_{i+1}}
\end{align*}
is a representative of the dense $B$-orbit in $\u_\alpha$.  In case
$\alpha = \alpha_n$, a representative of the dense $B$-orbit in
$\u_\alpha$ is
\[
e =  \left(\sum_{j = 1}^{n-1} e_{j,j+1} - e_{-(j+1),-j}\right) +
(e_{n-1,-n} - e_{n,-n-1})
= \left(\sum_{j = 1}^{n-1} e_{\alpha_j}\right) +
e_{\alpha_{n-1}+2\alpha_n}.
\]

\smallskip

Now consider $G = \Sp_{2n}(k)$.  For $\alpha = \alpha_1$, we take
\[
e = \left(\sum_{j = 2}^{n-1} e_{j,j+1} - e_{-(j+1),-j}\right) +
e_{n,-n} + e_{1,-1} =  \left(\sum_{j = 2}^n e_{\alpha_j}\right) +
e_\rho;
\]
and for $\alpha = \alpha_n$ we take
\[
e = \left(\sum_{j = 1}^{n-1} e_{j,j+1} - e_{-(j+1),-j}\right) +
e_{n-1,-(n-1)} = \left(\sum_{j = 1}^{n-1} e_{\alpha_j}\right) +
e_{2\alpha_{n-1} + \alpha_n}.
\]

\smallskip

Finally, we consider the case $G = \SO_{2n}(k)$.  For $\alpha =
\alpha_1$,
\begin{align*}
e = & \left(\sum_{j = 2}^{n-1} e_{j,j+1} - e_{-(j+1),-j}\right) +
(e_{n-1,-n} - e_{n,-(n-1)}) + (e_{1,-n} - e_{n,-1}) \\
= & \left(\sum_{j = 2}^n e_{\alpha_j}\right) + e_{(\alpha_1 +
\alpha_2 + \dots + \alpha_{n-2}) + \alpha_n}
\end{align*}
is a representative of the dense $B$-orbit.  We finish with the case
$\alpha = \alpha_{n-1}$; the case $\alpha = \alpha_n$ is equivalent, by
symmetry.  A representative of the dense $B$-orbit in $\u_\alpha$ is
\begin{align*}
e = & \left(\sum_{j = 1}^{n-2} e_{j,j+1} - e_{-(j+1),-j}\right) +
(e_{n-1,-n} - e_{n,-(n-1)}) + (e_{n-2,-(n-1)} - e_{n-1,-(n-2)}) \\
= & \left(\sum_{j = 1}^{n-2} e_{\alpha_j}\right) + e_{\alpha_n} +
e_{\alpha_{n-2} + \alpha_{n-1} + \alpha_n}.
\end{align*}
This completes the proof of Theorem \ref{T:richardson} and so of
Theorem \ref{T:fibre}.

\section{Subregular orbital varieties}
\label{S:proof:orbital}

In this section we apply Theorem \ref{T:richardson} to classify all
subregular orbital varieties admitting a finite number of
$B$-orbits. Let $\cO$ be the subregular nilpotent $G$-orbit in $\g$.
Since $\cO$ is the Richardson class of all semisimple rank $1$
parabolic subgroups $P_\alpha$ of $G$, it follows that the orbital
varieties of $\cO \cap \u$ coincide with the nilradicals $\u_\alpha$
of the $P_\alpha$ $(\alpha \in \Pi)$. Observe that Theorem
\ref{T:richardson} gives a classification of all instances when $B$
acts on $\u_\alpha$ with a dense orbit: for, this is equivalent to
$B$ admitting a dense orbit in $\cO \cap \u_\alpha$ which in turn is
equivalent to $B$ acting on $\cO \cap \u_\alpha$ with a finite
number of orbits, as explained at the end of the introduction.
However, in contrast to Problems  \ref{Prob:fibre} and
\ref{Prob:richardson} in the subregular case, the questions of
finiteness for the number of orbits versus that of the existence of
a dense orbit are not equivalent for the action of $B$ on the
orbital varieties $\u_\alpha$ in $\overline{\cO \cap \u}$.

\begin{thm}
\label{T:orbital}
There is a finite number of $B$-orbits in
$\u_\alpha$ if and only if
\begin{itemize}
\item[(i)]
$G$ is of type $A_r$ for $r \le 4$; or $G$ is of type $A_5$ and $\alpha \in
\{\alpha_1, \alpha_3, \alpha_5\}$;
\item[(ii)]
$G$ is of type $B_2$; or $G$ is of type $B_3$ and $\alpha =
\alpha_2$;
\item[(iii)]
$G$ is of type $C_2$; or $G$ is of type $C_3$ and $\alpha \in
\{\alpha_1, \alpha_3\}$;
\item[(iv)]
$G$ is of type $G_2$ and $\alpha = \alpha_2$.
\end{itemize}
\end{thm}

\begin{proof}
It was shown by Kashin in \cite{Ka} that $B$ acts on $\u$ with a
finite number of orbits if and only if $G$ is of type $A_r$ for $r
\le 4$ or $B_2$. So clearly, $B$ acts on each $\u_\alpha$ with a
finite number of orbits in these cases. This classification was
extended in \cite[Cor.\ 1.4]{PR} to the case of minimal parabolic
subgroups $P_\alpha$ of $G$. The results in \cite{PR} are stated and
proved assuming that $\Char k = 0$; they are also valid provided
$\Char k$ is a good prime for $G$, cf.\ \cite{Ro}. Accordingly,
$P_\alpha$ acts on $\u_\alpha$ with a finite number of orbits if and
only if $G$ is of type $A_r$ for $r \le 5$, $B_r$ for $r \le 3$,
$C_r$ for $r \le 3$, $D_4$, or $G_2$. Obviously, if $P_\alpha$
already acts on $\u_\alpha$ with infinitely many orbits, so does
$B$. Clearly, if $B$ acts on $\u_\alpha$ with a finite number of
orbits, then $B$ has only finitely many orbits in $\cO \cap
\u_\alpha$. We thus infer from Theorem \ref{T:richardson}(a) and the
classification results from \cite{Ka} and \cite{PR} that, given $B$
acts on all of $\u$ with infinitely many orbits, the only cases that
need consideration for $B$ acting on $\u_\alpha$ are as follows:
each simple root for $A_5$ and $B_3$; $\alpha_1$ and $\alpha_3$ for
$C_3$; each end-node simple root for $D_4$, and $\alpha_2$ for
$G_2$.

In order to show that $B$ acts with an infinite number of orbits in
several of these cases, we employ a method already used in \cite{Ka}
and \cite{PR} which we recall now for convenience. Let $\n$ be a
$B$-submodule of $\u$, i.e.\ an ideal of $\b$ in $\u$. Let $N$ be
the connected unipotent normal subgroup of $B$ with Lie algebra
$\n$. The action of $B$ on $\n$ induces an action of $B$ on the
quotient $\n/[\n,\n]$ of $\n$ by its commutator subalgebra, and this
latter action of $B$ factors through $B/N$. Thus if $\dim B/N < \dim
\n/[\n,\n]$, then $B$ acts on $\n$ with an infinite number of
orbits. The idea in all the cases we consider below is to exhibit a
suitable $B$-submodule $\n$ of $\u$ which satisfies this inequality.
We collect the relevant information in Table \ref{t:1} below where
we list $\n$ by means of the simple roots $\beta$ so that $\n$ is
generated by the root spaces $\g_\beta$ as a $B$-module. We omit the
details.

\begin{table}[h!]
\renewcommand{\arraystretch}{1.2}
\begin{tabular}{|cccc|}
\hline
Type of $G$ & $\n$ & $\dim B/N$ & $\dim \n/[\n,\n]$\\
\hline
$A_5$ & $\alpha_1, \alpha_3, \alpha_5$ & $7$ & $8$ \\
$B_3$ & $\alpha_2$ & $5$ & $6$ \\
$C_3$ & $\alpha_1, \alpha_3$ & $4$ & $5$ \\
$D_4$ & $\alpha_2$ & $7$ & $8$ \\
$G_2$ & $\alpha_2$ & $3$ & $4$ \\
\hline
\end{tabular}
\vspace{2mm}
\caption{Some ideals $\n$ in $\b$} \label{t:1}
\end{table}
Note that if $\alpha$ is a simple root of $G$ not among the roots of
the generating root spaces for $\n$ from Table \ref{t:1}, then $\n
\subseteq \u_\alpha$. It follows from the discussion above and the
information given in Table \ref{t:1} that then $B$ acts on
$\u_\alpha$ with an infinite number of orbits. This immediately
rules out
the remaining cases for $D_4$, as well as the cases $\alpha_2,
\alpha_4$ for $A_5$ and $\alpha_1, \alpha_3$ for $B_3$.
Consequently, it remains to be checked that each of the cases for
$A_5$, $B_3$, $C_3$ and $G_2$ listed in the statement is
indeed an instance when $B$ acts on $\u_\alpha$ with a finite number
of orbits. It turns out that for $G$ of type $A_5$, $B_3$, $C_3$,
$G_2$ there is a unique one-parameter family of $B$-orbits in $\u$
and moreover that this family is dense precisely in the ideal $\n$
given in Table \ref{t:1}; see Tables 1 and 2 in \cite{BH}, see also
the proof of \cite[Prop.\ 4.8]{HR} for type $A_5$. It thus follows
that $B$ acts with a finite number of orbits in $\u_\alpha$ for the
remaining instances listed in the statement, as then
$\u_\alpha$ does not meet this infinite family of orbits in $\n$.
\end{proof}

\begin{rem}
\label{rem:numbers}
In each of the finite instances of Theorem \ref{T:orbital}
the number of $B$-orbits in $\u_\alpha$ can be determined from the
computation of a list of representatives of
the $B$-orbits in all of $\u$ by means of
the algorithm from \cite{BH}.  Apart from the $A_5$ case,
the list of $B$-orbits
is given in \cite[Table 2]{BH}; the list for $G$ of type $A_5$ was made
available to us by W.~Hesselink.
If $G$ is of type $B_3$, then there are $23$ $B$-orbits in
$\u_{\alpha_2}$, if $G$ is of type $C_3$, there are $24$ $B$-orbits
in $\u_{\alpha_1}$ and $21$ $B$-orbits in $\u_{\alpha_3}$, while if
$G$ is of type $G_2$, there are $8$ $B$-orbits in $\u_{\alpha_2}$.
For $G$ of type $A_5$, there are $185$ $B$-orbits on $\u_{\alpha_1}$
and $200$ $B$-orbits on $\u_{\alpha_3}$.
\end{rem}

\bigskip

{\bf Acknowledgments}:
This research was funded in part by EPSRC
grant EP/D502381/1. Part of the research for this paper was carried
out while the authors were staying at the Mathematical Research
Institute Oberwolfach supported by the ``Research in Pairs''
programme. The first author would like to thank New
College, Oxford for financial support whilst the research was
carried out.
We are grateful to W.\ Hesselink for making available the
computations of the $B$-orbits in the $A_5$ cases that were used
in Remark \ref{rem:numbers}.


\bigskip

\begin{thebibliography}{00}

\bibitem{BR}
P.~Bardsley, R.~W.~Richardson, \emph{\'Etale slices for algebraic
transformation groups in characteristic $p$}, Proc.\ London Math.\
Soc.\ (3) \textbf{51} (1985), no.\ 2, 295--317.

\bibitem{Bo}
A.~Borel,
\emph{Linear algebraic groups}, Graduate Texts in Mathematics {\bf
126}, Springer-Verlag 1991.

\bibitem{BB}
W.~Borho, J.-L.~Brylinski,
\emph{Differential operators on
homogeneous spaces}, Invent.\ Math.\ {\bf 80} (1985), 1--68.

\bibitem{Bou}
N.~Bourbaki,
\emph{Groupes et alg\`{e}bres de Lie}, Chapitres 4, 5
et 6, Hermann, Paris, 1975.

\bibitem{Br}
M.~Brion,
\emph{Quelques propri\'et\'es des espaces homog\'enes
sph\'eriques}, Man.\ Math., \textbf{99} (1986), 191--198.

\bibitem{BH}
H.~B\"urgstein, W.~H.~Hesselink, \emph{Algorithmic orbit
classification for some Borel group actions}, \\ Comp.~Math.
\textbf{61} (1987), 3--41.

\bibitem{Ca}
R.~W.~Carter,
\emph{Finite groups of Lie type. Conjugacy classes and
complex characters}, Pure and Applied Mathematics, New York, 1985.

\bibitem{DLP}
C.~De Concini, G.~Lusztig, C.~Procesi, \emph{Homology of the
zero-set of a nilpotent vector field on a flag manifold}, J. Amer.
Math. Soc. \textbf{1} (1988), no.\ 1, 15--34.

\bibitem{FR}
R.~Fowler, G.~R\"ohrle,
\emph{Spherical nilpotent orbits in good characteristic}, in preparation.

\bibitem{Go}
S.~M.~Goodwin,
\emph{Algorithmic testing for dense orbits of Borel subgroups}, J.\
Pure Appl.\ Algebra, \textbf{197} (2005), no.\ 1--3, 171--181.

\bibitem{Go2}
\bysame, 
\emph{ On the conjugacy classes in maximal unipotent subgroups of
simple algebraic groups}, Transform.\ Groups, {\bf 11} (2006), no.\
1, 51--76.

\bibitem{GH}
S.~M.~Goodwin, L.~Hille,
\emph{Prehomogeneous spaces for Borel subgroups of general linear
groups}, Transform.\ Groups, to appear.

\bibitem{He}
W.~H.~Hesselink,
\emph{Singularities in the nilpotent scheme of a
classical group}, Trans. Amer. Math. Soc. \textbf{222} (1976),
1--32.

\bibitem{HR}
L.~Hille, G.~R\"ohrle, \emph{A classification of parabolic subgroups
of classical groups with a finite number of orbits on the unipotent
radical}, Transform.\ Groups \textbf{4} (1998), no.\ 1, 317--337.

\bibitem{HS}
R.~Hotta, N.~Shimomura,
\emph{The fixed-point subvarieties of
unipotent transformations on generalized flag varieties and the
Green functions},
Math. Ann. \textbf{241} (1979), no.\ 3, 193--208.

\bibitem{Ja}
J.~C.~Jantzen, \emph{Nilpotent orbits in representation theory}, Lie
theory, 1--211, Progr. Math., \textbf{228}, Birkh\"auser, Boston,
MA, 2004.

\bibitem{Jo}
A.~Joseph,
\emph{On the variety of a highest weight module}, J.\
Algebra {\bf 88} (1984), 238--278.

\bibitem{Ka}
V.V.~Kashin,
\emph{Orbits of adjoint and coadjoint actions of Borel
subgroups of semisimple algebraic groups},  (in Russian) Problems in
Group Theory and Homological algebra, Yaroslavl' (1990), 141--159.

\bibitem{Kn}
F.~Knop, \emph{On the set of orbits for a Borel subgroup}, Comment.\
Math.\ Helv.\ \textbf{70} (1995), 285--309.

\bibitem{Lu}
G.~Lusztig,
\emph{Green polynomials and singularities of unipotent classes},
Adv.\ in Math.\ {\bf 42} (1981), no.\ 2, 169--178.

\bibitem{Pa}
D.~Panyushev, \emph{Complexity and nilpotent orbits}, Manuscripta
Math.\ {\bf 83} (1994), no.\ 3--4, 223--237.

\bibitem{PR}
V.~Popov, G.~R\"ohrle,
\emph{On the number of orbits of a
parabolic subgroup on its unipotent radical}, in ``Algebraic Groups
and Lie Groups'' Australian Mathematical Society Lecture Series, ed.
G.I.~Lehrer, \textbf{9}, (1997), 297--320.

\bibitem{Ri}
R.W.~Richardson,
\emph{Conjugacy classes in parabolic subgroups of
semisimple algebraic groups}, Bull. London Math. Soc. \textbf{6}
(1974), 21--24.

\bibitem{Ro}
G.~R\"ohrle, \emph{On the modality of parabolic subgroups of linear
algebraic groups}, Manuscripta Math.\ \textbf{98} (1999), 9--20.

\bibitem{Sl}
P.~Slodowy,
\emph{Simple singularities and simple algebraic groups},
Lecture Notes in Mathematics, \textbf{815} Springer-Verlag,
Berlin-New York, 1980.

\bibitem{Sp}
N.~Spaltenstein,
\emph{Classes unipotentes et sous-groupes de
Borel}, Lecture Notes in Mathematics, \textbf{946} Springer-Verlag,
Berlin-New York, 1982.

\bibitem{Spr}
T.~A.~Springer,
\emph{The unipotent variety of a semi-simple group}, Algebraic
Geometry (Internat.\ Colloq., Tata Inst.\ Fund.\ Res., Bombay,
1968), Oxford Univ.\ Press, London, 1969, pp.\ 373--391.

\bibitem{Spr2}
\bysame, 
\emph{Trigonometric sums, Green functions of finite groups and
representations of Weyl groups}, Invent.\ Math.\ {\bf 36} (1976),
173--207.

\bibitem{SS}
T.~A.~Springer, R.~Steinberg,
\emph{Conjugacy classes},
Seminar on Algebraic Groups and Related Finite Groups,
Lecture Notes in Mathematics, \textbf{131}, (1968-69), 167--266,
Springer-Verlag, Berlin-New York, 1970.

\bibitem{St}
R.~Steinberg, \emph{Conjugacy classes in algebraic groups}, Lecture
Notes in Mathematics, \textbf{366}, Springer-Verlag, Berlin-New
York, 1974.

\bibitem{St2}
\bysame, 
\emph{On the desingularization of the unipotent variety}, Invent.\
Math.\ {\bf 36} (1976), 209--224.

\bibitem{Vi}
E.~B.~Vinberg, \emph{Complexity of actions of reductive groups},
Funct.\ Anal.\ Appl.\ \textbf{20} (1986), 1--11.

\end{thebibliography}
\end{document}